\documentclass{amsart} 

\newtheorem{theorem}{Theorem}

\newtheorem{proposition}[theorem]{Proposition}
\newtheorem{corollary}[theorem]{Corollary}

\theoremstyle{definition}

\newcommand{\Z}{\mathbb{Z}}
\newcommand{\R}{\mathbb{R}}

\newcommand{\T}{\mathbb{T}}

\newcommand{\lw}[1]{\smash{\lower6pt\hbox{#1}}} 

\usepackage{graphics} 

\begin{document}


\title[Ribbon-moves and Khovanov-Jacobsson numbers]
{Ribbon-moves for 2-knots with 1-handles attached\\
and Khovanov-Jacobsson numbers} 


\author{J. Scott Carter} 
\address{Department of Mathematics, 
University of South Alabama, Mobile, AL 36688, U.S.A.} 
\email{carter@jaguar1.usouthal.edu} 

\author{Masahico Saito} 
\address{Department of Mathematics, 
University of South Florida, Tampa, FL 33620, U.S.A.} 
\email{saito@math.usf.edu} 

\author{Shin Satoh} 
\address{Graduate School of Science and Technology, 
Chiba University, Yayoi-cho 1-33, 
Inage-ku, Chiba, 263-8522, Japan 
(Department of Mathematics, University of South Florida, 
April 2003--March 2005)}
\email{satoh@math.s.chiba-u.ac.jp}

\renewcommand{\thefootnote}{\fnsymbol{footnote}}
\footnote[0]{2000 {\it Mathematics Subject Classification}. 
Primary 57Q45; Secondary 57Q35.}

\keywords{Khovanov homology, 2-knot, 
ribbon-move, twist-spun knot, crossing change.}

\maketitle

\begin{abstract} 
We prove that a crossing change along a double point circle 
on a $2$-knot 
is realized by ribbon-moves for
a knotted torus obtained 
from the $2$-knot by attaching  a $1$-handle. 
It follows that any $2$-knots for which 
the crossing change is an unknotting operation, 
such as ribbon $2$-knots and twist-spun knots, 
have  trivial 
Khovanov-Jacobsson number. 
\end{abstract}

\begin{figure}[htb]
\begin{center}
\includegraphics{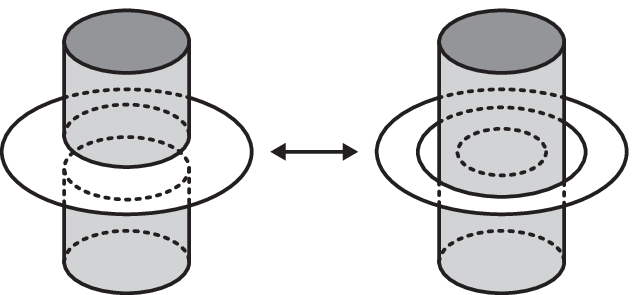}
\caption{}
\label{fig01}
\end{center}
\end{figure}

A {\it surface-knot} or {\it -link} 
is a closed surface 
embedded in $4$-space $\R^4$ locally flatly. 
Throughout this note, 
we always assume that all surface-knots are oriented. 
A {\it ribbon-move} (cf.~\cite{Og}) is a local operation 
for (a diagram of) a surface-knot 
as shown in Figure~\ref{fig01}. 
We say that surface-knots $F$ and $F'$ are 
{\it ribbon-move equivalent}, 
denoted by $F\sim F'$, 
if $F'$ is obtained from $F$ 
by a finite sequence of ribbon-moves.

\begin{figure}[htb]
\begin{center}
\includegraphics{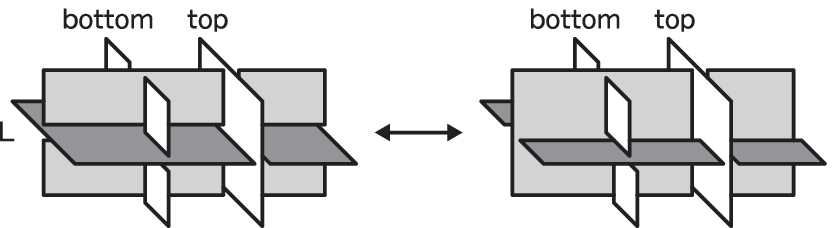}
\caption{}
\label{fig02}
\end{center}
\end{figure}

The ribbon-move is a special case of the {\it crossing change}: 
Assume that a surface-knot $F$ has 
a double point circle $L$ in a diagram 
such that 
(i) $L$ has no self-intersection, and 
(ii) at every triple point on $L$, 
the  sheet transverse to $L$ is either
 top or bottom (not middle). 
 The condition (i) means that $L$ 
 does not go through the same triple point twice.
When $L$ satisfies these conditions, we can perform 
a crossing change along $L$ 
by exchanging the roles 
of 
over- and under-sheets as indicated in Figure~\ref{fig02} 
(cf. \cite{Yas}). 
See \cite{CSbook} for details on diagrams of surface-knots.

For a $2$-knot $K$ (a knotted sphere in $\R^4$), 
a crossing change is not necessarily realized 
by ribbon-moves; 
indeed, a ribbon-move does not change 
the Farber-Levine pairing of $K$ 
but a crossing change might (cf. \cite{Og}). 
On the other hand, 
when we consider the $\T^2$-knot (knotted torus in $\R^4$) $K+h$ 
obtained from $K$ by attaching  a $1$-handle $h$ on $K$, 
we obtain the following. 

\begin{theorem}\label{thm1} 
Let $K$ and $K'$ be $2$-knots such that 
$K'$ is obtained from $K$ by a crossing change. 
Then for any $1$-handles $h$ and $h'$ on $K$ and $K'$, 
respectively, 
the $\T^2$-knot $K+h$ is ribbon-move equivalent to $K'+h'$. 
\end{theorem} 

\begin{figure}[htb]
\begin{center}
\includegraphics{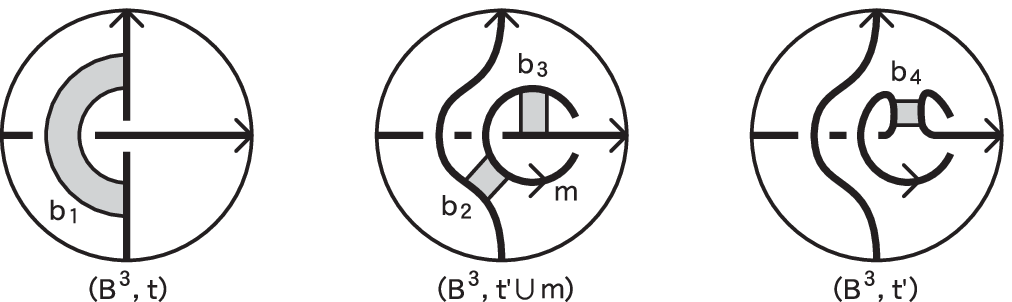}
\caption{}
\label{fig03}
\end{center}
\end{figure}

\begin{proof} 
Along the double point circle $L$ for which 
we perform the crossing change, there is 
  a neighborhood $N$ identified with 
$(B^3,t)\times S^1$, 
where $(B^3,t)$ is a tangle 
with two strings as shown in the left of Figure~\ref{fig03}. 
In the figure, the orientations of tangles are induced from that of $K$,
and all bands are attached in an orientation-compatible manner.
For an interval $I$ in $S^1$, 
we take a $1$-handle $h_1=b_1\times I$ on $K$, 
where $b_1$ is a band as indicated in the figure. 

We observe that $K+h_1$ is ambient isotopic to 
$(K'\cup T)+h_2$ (cf.~\cite{Sa2}), 
where $T=m \times S^1$ is a 
$\T^2$-knot 
linking with $K'$, 
and the $1$-handle $h_2=b_2\times I$ 
connects between $K'$ and $T$. 
See the center of Figure~\ref{fig03}. 

Consider a $1$-handle $h_3=b_3\times I$ 
on $K'\cup T$. 
Since both of $h_2$ and $h_3$ 
connect between $K'$ and $T$, 
the $\T^2$-knot $(K'\cup T)+h_2$ 
is ribbon-move equivalent to $(K'\cup T)+h_3$. 

Finally we see that 
$(K'\cup T)+h_3$ is 
ambient isotopic to $K'+h_4$, 
where $h_4=b_4\times I$ is the $1$-handle on $K'$ 
as shown in the right of the figure. 
Thus we obtain 
$$K+h\sim K+h_1=
(K'\cup T)+h_2\sim (K'\cup T)+h_3
=K'+h_4\sim K'+h'.$$
This completes the proof. 
\end{proof}

We say that the crossing change is an {\it unknotting operation} 
for a surface-knot $F$ 
if the trivial surface-knot  is obtained from $F$ 
by a finite sequence of crossing changes. 
It is still unknown
whether the crossing change is an unknotting operation 
for {\it any} surface-knot. 

Khovanov \cite{Kh} introduced a categorification of
the Jones polynomial, 
that is, a chain complex for a given classical knot diagram
such that its graded Euler characteristic
is the Jones polynomial. 
Khovanov \cite{Kh4}  and Jacobsson \cite{Ja} proved 
that it defines an invariant for cobordisms (relative to boundary diagrams).
Specifically, a cobordism between two knot diagrams gives rise to
a chain map (we call it a Khovanov-Jacobsson homomorphism) 
between corresponding chain complexes, 
that is invariant 
under equivalence of cobordisms of diagrams. 
See also \cite{Dror}.
In particular,  
a diagram of a $\T^2$-knot 
is a cobordism between empty diagrams, 
giving rise to a homomorphism $\Z \rightarrow \Z$ 
 defined up to sign,  a multiplication by a constant.
We call this constant the {\it Khovanov-Jacobsson number}.

\begin{theorem}\label{thm2} 
Let $K$ be a $2$-knot for which 
the crossing change is an unknotting operation. 
Then for any $1$-handle $h$ on $K$, 
the $\T^2$-knot $K+h$ has the trivial Khovanov-Jacobsson number. 
\end{theorem} 

\begin{figure}[htb]
\begin{center}
\includegraphics{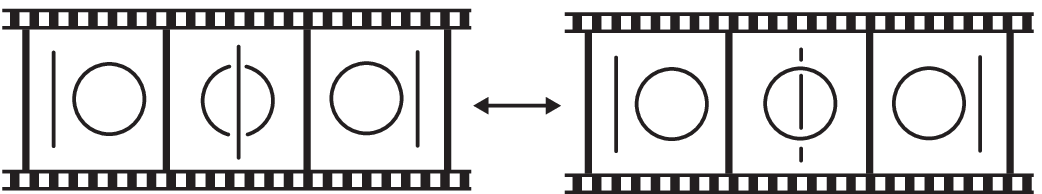}
\caption{}
\label{fig04}
\end{center}
\end{figure}

\begin{proof} 
Let $K_0$ be the trivial $2$-knot 
and $h_0$ the trivial $1$-handle on $K_0$. 
By assumption and Theorem~\ref{thm1}, the $\T^2$-knot 
$K+h$ is ribbon-move equivalent to 
$K_0+h_0$, 
which is the trivial $\T^2$-knot. 

Consider two movies as shown in Figure~\ref{fig04}. 
It is seen from the definitions \cite{Dror,Ja}  that 
the corresponding Khovanov-Jacobsson homomorphisms 
$H^*(|\bigcirc)\rightarrow H^*(\bigcirc |)$ 
are the same for these movies. 
This implies that 
a ribbon-move does not change 
the Khovanov-Jacobsson number. 
Hence the $\T^2$-knot $K+h$ has 
the same number as that of 
the trivial $\T^2$-knot $K_0+h_0$. 
\end{proof} 

By Theorem~\ref{thm2}, 
if there is a $2$-knot $K$ 
such that the Khovanov-Jacobsson number of $K+h$ is non-trivial, 
then the crossing change  is not 
an unknotting operation for $K$. 
However, 
we have no such examples at present. 

\begin{corollary}\label{cor3} 
Let $K$ be a ribbon $2$-knot or twist-spun knot. 
Then for any $1$-handle $h$ on $K$, 
the $\T^2$-knot $K+h$ 
has 
trivial Khovanov-Jacobsson number. 
\end{corollary} 

\begin{proof} 
This follows from Theorem~\ref{thm2} and the fact that 
the crossing change is an unknotting operation for 
every ribbon $2$-knot or twist-spun knot 
(cf. \cite{AS, Sa1}). 
\end{proof} 

We say that a surface-knot is {\it pseudo-ribbon} \cite{Ka} 
if it has a diagram without triple points. 
The notions of ribbon and pseudo-ribbon $2$-knots 
are the same \cite{Ya} (see also \cite{KS}). 
On the other hand, 
for $\T^2$-knots, 
they are not coincident in the sense that 
the family of pseudo-ribbon $\T^2$-knots 
properly contains that of ribbon $\T^2$-knots.

\begin{proposition}\label{prop4} 
Any pseudo-ribbon $\T^2$-knot 
has trivial Khovanov-Jacobsson number. 
\end{proposition} 

\begin{proof} 
By the results of Teragaito \cite{Te} and Shima \cite{Sh}, 
every pseudo-ribbon $\T^2$-knot $T$ 
is (i) a ribbon $\T^2$-knot, or 
(ii) a $\T^2$-knot obtained from 
a split union of a Boyle's turned $\T^2$-knot $T'$ \cite{Bo} 
and a trivial $2$-link $U=U_1\cup U_2\cup\dots\cup U_n$ 
by surgery along $1$-handles 
$h_1,h_2,\dots,h_n$ for some $n\geq 0$, 
where each $h_i$ connects between $T'$ and 
$h_i$ $(i=1,2,\dots,n)$. 

For the case (i), 
there is a ribbon $2$-knot $K$ and a $1$-handle $h$ 
such that $T=K+h$. 
Hence the conclusion follows from Corollary~\ref{cor3}. 

For the case (ii), 
we see that $T=(T'\cup U)+(\bigcup_{i=1}^n h_i)$ 
is ribbon-move equivalent to $T'$. 
We consider two movies for a classical knot diagram $D$ 
in a plane, 
one of which keep $D$ still 
and the other twists $D$ 
by a $2\pi$-rotation of the plane. 
Then it follows from the definitions \cite{Dror, Ja, Kh4}
that the corresponding Khovanov-Jacobsson homomorphisms 
$H^*(D)\rightarrow H^*(D)$ 
are the same for these movies. 
This implies that $T'$ has the same 
Khovanov-Jacobsson number 
as that of a non-turned (that is, just spun) $\T^2$-knot, 
which is ribbon. 
Hence this case reduces to (i). 
\end{proof}

\section*{Achknowledgments}

The first, second, and third authors are 
partially supported by 
NSF Grant DMS $\#0301095$, 
NSF Grant DMS $\#0301089$, 
and JSPS Postdoctoral Fellowships for Research Abroad, 
respectively. 
The third author expresses his gratitude 
for the hospitality of the University of South Florida 
and the University of South Alabama.


\end{document}